\magnification=\magstep1 
\baselineskip=12pt 

\def\ind {{\lim \limits_{\longrightarrow}}} 
\def \pro {{\lim\limits_{\longleftarrow}}} 
\input amstex 
\documentstyle{amsppt} 
\overfullrule=0pt 
\vsize=500pt 

\document

\def\gen{{\frak{g}}}

\def\lra{\,\longrightarrow\,}

\def\l{{\lambda}}

\def\eps{{\varepsilon}}

\def\1b{{\bold 1}}

\def\L{{\roman L}}

\def\ft{{\text{ft}}}

\def\gr{\roman{gr}} 
 
\def\Spec{\roman{Spec\,}} 
\def\Spf{\roman{Spf\,}}

\def\red{{\text{red}}} 
 
\def\ft{{\text{ft}}} 
\def\sqr{{\!\sqrt{\,\,\,}}} 
 
\def\ch{{\text{ch}}}

\def\Schb{{\bold{Sch}}}

\def\Ischb{{\bold{Isch}}}

\def\Aut{\roman{Aut}} 
\def\Ext{\text{Ext}}

\def\Hom{\roman{Hom}} 

\def\Der{\roman{Der}}

\def\Res{\text{Res}}

\def\End{\text{End}\,} 
 
\def\tr{\text{tr}\,}

\def\Ker{\text{Ker}\,}

\def\Im{\text{Im}\,}

\def\AA{{\Bbb A}} 
 
\def\CC{{\Bbb C}}

\def\RR{{\Bbb R}}

\def\ZZ{{\Bbb Z}}

\def\Ac{{\Cal A}} 
 
\def\Cc{{\Cal C}} 
\def\Dc{{\Cal D}} 
\def\Ec{{\Cal E}} 
\def\Fc{{\Cal F}}

\def\Kc{{\Cal K}} 
\def\Lc{{\Cal L}}

\def\Oc{{\Cal O}}

\def\Rc{{\Cal R}}

\def\Vc{{\Cal V}} 
\def\Wc{{\Cal W}}

\def\and{{\quad\text{and}\quad}}

\def\qed{\hfill $\sqcap \hskip-6.5pt \sqcup$}

\def\simto{\,{\buildrel\sim\over\to}\,} 
\def\lra{{{\longrightarrow}}}

\def\u1{{\underline 1}}

\centerline{\bf FORMAL LOOPS III :}
\centerline{\bf FACTORIZING FUNCTIONS AND THE RADON TRANSFORM} 

\vskip .5cm 

\centerline {\bf M. Kapranov, E. Vasserot} 

\vskip 1cm

\vskip 1cm
\centerline{\bf 0. Introduction.} 

\vskip3mm

\noindent{\bf (0.1)}
Let $X$ be  a $C^\infty$-manifold and
$\eta$ be a smooth 1-form on $X$. Its Radon
transform is the function $\tau(\eta)$ on the free loop space
$L(X) = C^\infty (S^1, X)$ whose value at  $\gamma: S^1\to X$
is $\int_{S^1} \gamma^*\eta$. This ``universal'' setting
includes many classical instances of integral transforms over families of
curves
such as the original Radon transform involving integration
over straight lines in $\RR^3$.

The problem of describing the range of the transform $\tau$
was studied by  J.-L. Brylinski [B] who, in  the case $X=\RR^N$,
characterized it by
 a system of differential equations
similar to the
A-hypergeometric system [GKZ]. One part of these equations
expresses the following obvious property:

\vskip2mm

\noindent{(0.1.1)} 
$\tau(\eta)$ is invariant under reparametrizations of $S^1$.

\vskip2mm

Suppose that $X$ is a complex manifold and
$\eta$ is holomorphic. Then $\tau(\eta)$ satisfies  two more 
properties:

\vskip2mm

\noindent(0.1.2) 
$\tau(\eta)=0$ on $L^0(X)$, the
space  of loops extending holomorphically into the unit disk.

\vskip1mm

\noindent(0.1.3) 
$\tau(\eta)$ is additive in holomorphic pairs of pants
(holomorphic maps of sphere minus 3 disks into $X$). Namely, the
value on $\tau(\eta)$ on the ``waist'' circle is equal to the sum
of values on the ``leg'' circles. 

\vskip .3cm

\noindent {\bf (0.2)} 
In the present paper we further assume that $X$ is a smooth
algebraic variety over $\CC$ and relate
 the Radon transform  to the theory of vertex
operator algebras and of related objects, factorization algebras
introduced by A. Beilinson and V. Drinfeld [BD]. We
 replace $L^0(X)$ and $L(X)$
by the scheme $\Lc^0(X)$ of formal arcs [DL] and the
ind-scheme $\Lc(X)$ of formal loops [KV1]. Informally, they are obtained by
replacing the unit disk and the circle by $\Spec\, \CC[[t]]$ and
$\Spec \, \CC((t))$. Then for a regular
1-form $\eta$ we have a function $\tau(\eta)$ on $\Lc(X)$ vanishing
on $\Lc^0(X)$. In fact, for any closed 2-form $\omega$
the Radon transform $\tau (d^{-1}(\omega))$ of any analytic
primitive $d^{-1}(\omega)$  of $\omega$ is defined algebraically and depends
only on $\omega$. If $\omega$ is nondegenerate, then $\tau(d^{-1}(\omega))$
is a version of the symplectic action functional on the space of loops. 

In [KV1] we showed how the space $\Lc(X)$ provides a geometric
construction of a particular sheaf of vertex algebras on $X$,
the chiral de Rham complex $\Omega^{\ch}_X$ from [MSV]. It was realized as a
semiinfinite de Rham complex of a particular D-module on $\Lc(X)$
``consisting'' of distributions  supported on $\Lc^0(X)$. In particular,
the sheaf $\Oc_{\Lc(X)}$ of functions on $\Lc(X)$ acts on $\Omega^{\ch}_X$
by multiplication.
 Our main result is the following.

\proclaim{(0.2.1) Theorem}Let $f$ be a function on $\Lc(X)$ vanishing on
$\Lc^0(X)$. Then the following are equivalent:

(i) $f$ has the form $\tau(d^{-1}(\omega))$ for a closed 2-form
$\omega$.

(ii) The operator of  multiplication by $\exp(f)$  in
$\Omega^{\ch}_X$  is an automorphism of vertex algebras. 

(iii) The function $f$ is 2-factorizing (see below).
\endproclaim

\vskip .3cm

\noindent {\bf (0.3)} The condition for $f$ to be 2-factorizing
can be seen as an infinitesimal version of the conditions (0.1.1) 
and (0.1.3).
First, we require that $f$ is invariant under reparametrizations
of $\Spf\,\CC[[t]]$, see (1.2.8). 
If $C$ is a smooth curve and $c\in C$ a point, then
there is a version $\Lc_{C,c}(X)$ of $\Lc(X)$ involving the punctured
formal disk on $C$ near $c$. Reparametrization invariance of $f$
implies
 that it gives rise to a function $f_{C,c}$ on $\Lc_{C,c}(X)$. 

Now, in [KV1] we introduced the total space $\Lc_C(X)\buildrel\pi\over
\longrightarrow
C$ of all the $\Lc_{C,c}(X)$ and showed that it has a structure of
a factorization monoid. This structure includes, in particular,
a family $\Lc_{C^2}(X)\to C^2$ whose fiber over $(c_1, c_2)$ with
$C_1\neq c_2$ is $\Lc_{C, c_1}(X)\times\Lc_{C, c_2}(X)$
while the fiber over $(c,c)$ is $\Lc_{C,c}(X)$. When $c_1, c_2$ merge
together to a point $c$, the punctured formal disks around $c, c_1$
and $c_2$ form an ``infinitesimal pair of pants''. A function
$f$ is called 2-factorizing if there exists a function
$f_{C^2}\in \Oc(\Lc_{C^2}(X))$ whose restriction on
$\pi^{-1}(c_1, c_2)$ with $c_1\neq c_2$ is $f_{C, c_1}+f_{C, c_2}$
while the restriction to $\pi^{-1}(c,c)$ is $f_{C,c}$. 
So this is an infinitesimal version of (0.1.3). The general concept
of factorizing functions (Definition 1.5.2) involves similar conditions
for any number of points merging together and turns out to
be equivalent to the 2-factorizability. 

\vskip .3cm

\noindent {\bf (0.4)} This paper is a part of a series [KV1-3] devoted to
 interpreting,
in geometric terms involving $\Lc(X)$,  the gerbe of chiral differential
operators (CDO) on $X$, see [GMS1-2]. Objects of this gerbe
are sheaves of vertex algebras similar to $\Omega^{\ch}_X$ but without
the fermionic variables. For any such object $\Ac$ the sheaf
$\underline{\Aut}(\Ac)$ was found in [GMS1]
to be canonically isomorphic to $\Omega^{2, cl}_X$, the sheaf of closed
2-forms.  This is precisely the domain of definition
of the Radon transform $\tau\circ d^{-1}$ from (0.2). 
In fact, our proof of the equivalence of (i) and (ii)
in Theorem 0.2.1 is based on this identification. 

In a subsequent paper [KV3] we will use Theorem 0.2.1 to relate
the gerbe of CDO with an appropriate version of the determinantal
gerbe of the tangent bundle of $\Lc(X)$. In other words, we will show
that the anomaly inherent in construction of CDO is
related to the determinantal anomaly for the loop space by the Radon
transform. Sheaves of CDO on $X$ form a gerbe with lien $\Omega^{2, cl}_X$
while the determinantal gerbe of $\Lc(X)$ has the lien $\Oc^\times_{\Lc(X)}$
On the level of liens the relation between the two gerbes associates
to a 2-form  $\omega\in\Omega^{2, cl}_X$ the invertible function
$S(\omega) = \exp(\tau d^{-1}(\omega))\in\Oc^\times_{\Lc(X)}$.

\vskip .3cm

Research of M.K. was partially supported by an NSF grant.

\vfill\eject

\centerline{\bf 1. Formal loop spaces and transgressions.} 

\vskip 1cm

\subhead{\bf (1.1) Conventions}\endsubhead
All rings will be assumed commutative and containing the field
$\CC$ of complex numbers. For a ring $R$ we denote by $R((t))^\sqr$
the subring in the ring $R((t))$ of formal Laurent series formed
by series $\sum_{n\gg -\infty}^\infty a_n t^n$ with $a_n$ nilpotent
for $n<0$, see [KV1]. 
Let $\Schb$ be the category of schemes over $\CC$,
and $\Ischb$ the category of ind-schemes with countable indexing.
Let $\Schb^\ft\subset\Schb$ 
be the full subcategory of schemes of finite type. 

\subhead{\bf (1.2) Reminder on $\Lc X$}\endsubhead

Let $X$ be a scheme of finite type over $\CC$.
We have then the scheme $\Lc^0X$ of formal arcs,
and the ind-scheme $\Lc X$ of formal loops in $X$, see \cite{KV1}.
They represent the following functors
$$\Hom(\Spec R, \Lc^0 X)=\Hom(\Spec R[[t]], X), \leqno(1.2.1)$$
$$\Hom(\Spec R, \Lc X)=\Hom(\Spec R((t))^\sqr, X). \leqno(1.2.2)$$
There is a diagram
$$X{\buildrel p\over\leftarrow}\Lc^0 X{\buildrel i\over\rightarrow}\Lc X,
\leqno(1.2.3)$$
with $p$ affine and $i$ a closed embedding.
Moreover, $\Lc X$ is an inductive limit of schmes $\Lc^\eps X$ that are
nilpotent extensions of $\Lc^0X$, and so it can be described as a locally 
ringed space $(\Lc^0X,\Oc_{\Lc X})$ where 
$$\Oc_{\Lc X}=\pro_\eps\Oc_{\Lc^\eps X}.$$
Notice that $\Oc_{\Lc X}$ has a natural topology (of the projective limit).

\vskip3mm

\noindent{\bf (1.2.4) Example.} 
Let $X=\AA^N$ with coordinates $b^1,b^2,...b^N$. Then
$$\Lc^0X=\Spec\CC\bigl[b^i_n;i=1,...N, n\ge 0\bigr],$$
where the $b^i_n$'s can be thought of as the coefficients of $N$
indeterminate Taylor series
$$b^i(t)=\sum_{n\ge 0}b^i_nt^n.$$
Further,
$$\Oc_{\Lc\AA^N}=\Oc_{\Lc^0\AA^N}\bigl[\bigl[b^i_N,i=1,...N, n<0\bigr]\bigr]$$
is the sheaf formed by formal power series in infinitely many variables
$b^i_n$, $n<0$.
By definition, such a series is a formal sum of monomials, each involving only
finitely many $b^i_n$, and the coefficients at the monomials can be arbitrary.
In particular
$$\Gamma(\Lc^0\AA^N,\Oc_{\Lc\AA^N})=
\CC\bigl[b^i_n,n\ge 0\bigr]\bigl[\bigl[b^i_n,n<0\bigr]\bigr].$$

\vskip1mm

Let
$$\Oc_{\Lc X|\Lc^0X}=\Ker\{\Oc_{\Lc X}\to\Oc_{\Lc^0X}\}.\leqno(1.2.5)$$
This is a sheaf of ideals in $\Oc_{\Lc X}$ consisting of topologically
 nilpotent
elements. Let also
$$\Oc^\times_{\Lc X|\Lc^0X}=\Ker\{\Oc^\times_{\Lc X}\to\Oc^\times_{\Lc^0X}\}.
\leqno(1.2.6)$$
Because of topological nilpotency of
$\Oc_{\Lc X|\Lc^0X},$
the exponential series defines an isomorphism
$$\exp:\Oc_{\Lc X|\Lc^0X}\to\Oc^\times_{\Lc X|\Lc^0X}.\leqno(1.2.7)$$
Let $K$ be the group scheme $\Aut\CC[[t]]$.
Explicitly
$$K=\Spec\CC[a_1^{\pm 1},a_2,a_3,...],\leqno (1.2.8)$$
where the $a_i$ can be thought of as coefficients of an indeterminate formal 
change of the variable
$$t\mapsto a_1t+a_2t^2+...$$
The group scheme $K$ acts on $\Lc X$ and $\Lc^0X$ 
by the above changes of the variable. 
We also denote by $\gen = \Der \, \CC[[t]] =\CC[[t]] d/dt$ 
the Lie algebra of derivations of $\CC[[t]]$. This algebra
acts on $\Lc^0(X), \Lc(X)$ in a way compatible with the action of $K$.
So we will speak about the action of the Harish-Chandra pair $(\gen, K)$.

\vskip .3cm

\subhead{\bf (1.3) The transgression}\endsubhead
Let $R=\pro_\alpha R_\alpha$ be a topological algebra represented as a
filtering
projective limit of discrete algebras $R_\alpha$.
The ring of Laurent series with coeffisients in $R$ is defined by
$$R((t))=\pro_\alpha R_\alpha((t)).\leqno (1.3.1)$$
An element of $R((t))$ can be viewed as a formal series
$\sum_{n=-\infty}^{+\infty}a_nt^n$,
$a_n\in R$, possibly infinite in both directions,
but with the condition $\lim\limits_{n\to -\infty} a_n=0$.

Recall from \cite{KV2, sec.~6.2} the evaluation map
which is a morphism of ringed space
$$ev=ev_X:(\Lc^0X,\Oc_{\Lc X}((t)))\to(X,\Oc_X), \leqno (1.3.2)$$
whose underlying morphism of topological spaces is the projection
$p:\Lc^0X\to X$.

\vskip3mm

\noindent{\bf (1.3.3) Example.}
Let $X=\AA^N$ as in Example 1.2.4. Then the morphism
$$H^0(ev):H^0(\AA^N,\Oc)\to H^0\bigl(\Lc^0\AA^N,\Oc_{\Lc\AA^N}((t))\bigl)$$
sends
$b^i$ to $\sum_{n=-\infty}^{+\infty}b_n^it^n$.

\vskip3mm

For any commutative algebra $R$ (topological or not) we denote by
$\Omega^1(R)$ the module of K\"ahler differentials 
(defined without taking into account the topology).
We also define $\Omega^i(R)=\bigwedge^i_R\Omega^1(R)$
(the algebraic exterior power).
We have then the de Rham differential
$d:\Omega^i(R)\to\Omega^{i+1}(R)$.

For any topological algebra $R$ as above we have the residue homomorphism
$$\Res=\Res^p:\Omega^p(R((t)))\to\Omega^{p-1}(R).$$
To define it, we denote, for each 
$a(t)=\sum a_nt^n\in R((t))$,
$$d_Ra(t)=\sum d(a_n)t^n\in\Omega^1(R)((t)),\quad
a'(t)=\sum na_nt^{n-1}\in R((t)).$$
Then we put
$$\Res^p(a_0(t)da_1(t)...da_p(t))=
\bigl(a_0(t)\sum_{i=1}^pd_Ra_1(t)...d_Ra_{i-1}(t)
a'_i(t)d_Ra_{i+1}(t)...d_Ra_p(t)\bigr)_{-1},\leqno (1.3.4)$$
where the subscript $(-1)$ means the coefficient at $t^{-1}$ in a series
from $\Omega^{p-1}(R)((t))$.

\proclaim{(1.3.5) Proposition}
The map $\Res^p$ is well defined and satisfies the following properties :

(a) If $\omega\in\Omega^p(R[[t]])$, then $\Res^p(\omega)=0$.

(b) If $\xi\in\gen$ and $L_\xi$ is its action on $\Omega^p\, R((t))$,
then $\Res ^p(\L_\xi \omega)=0$ (invariance of residue).


(c) $\Res^p(d\omega)=d\Res^p(\omega)$.
\endproclaim

The proof is standard.
Let
$$\Omega^p_{\Lc X|\Lc^0X}=\Ker\{\Omega^p_{\Lc X}\to\Omega^p_{\Lc^0X}\}.$$
We now define a morphism of sheaves
$$\tau=\tau^p:\Omega^p_X\to p_*\Omega^{p-1}_{\Lc X|\Lc^0X}\leqno(1.3.6)$$
called the transgression (or the Radon transform).
It is defined as the composition of
$$\Res: p_*\Omega^p\bigl(\Oc_{\Lc X}((t))\bigr)\to 
p_*\Omega^{p-1}(\Oc_{\Lc X})=p_*\Omega^{p-1}_{\Lc X}$$
and the pullback with respect to the evaluation map
$$ev^*:\Omega^p_X=\Omega^p(\Oc_X)\to p_*\Omega^p\bigl(\Oc_{\Lc X}((t))\bigr).$$
Proposition 1.3.5 implies

\proclaim{(1.3.7) Proposition}
(a) For any local section $\omega$ of $\Omega^p_X$ the form $\tau(\omega)$
lies in $p_*\Omega^{p-1}_{\Lc X|\Lc^0X}$.

(b) The form $\tau(\omega)$ is $\gen$-invariant.

(c) $\tau(d\omega)=d\tau(\omega)$.
In particular, if $\omega\in\Omega^1_X$ is of the form $\omega=df$ for
a local regular function $f$, then $\tau(\omega)=0$.
\endproclaim

\noindent
{\bf (1.3.8) Example.}
Let
$$\omega=\sum_{i=1}^Nf_i(b^1,...,b^N)db^i,
\quad
f_i\in\CC[b^1,...,b^N]$$
be a 1-form regular in $\AA^N$.
Then $\tau(\omega)$ is the element of
$$\CC[b^i_{\ge 0}][[b^i_{<0}]]=H^0(\Lc^0\AA^N,\Oc_{\Lc\AA^N})$$
defined by
$$\tau(\omega)=\sum_{i=1}^N
\Res\biggl(f_i\biggl(\sum_{n=-\infty}^\infty b_n^1t^n,...
\sum_{n=-\infty}^\infty b_n^Nt^n\biggr)\cdot\sum_{n=-\infty}^\infty
 nb_n^it^{n-1}dt\biggr).$$

\vskip .3cm

\subhead{\bf (1.4) The map $\tau d^{-1}$ and the symplectic action functional}
\endsubhead

Let  $X$ be a smooth algebraic variety over $\CC$, and $\omega\in\Omega^p(X)$
 be a closed 
$p$-form on $X$. Then locally on the transcendental topology $\omega=d\eta$
for a $(p-1)$-form $\eta$ with analytic coefficients.
In this subsection we observe that $\tau^{p-1}(\eta)$ is in fact an
(algebraic) section of 
$p_*\Omega^{p-2}_{\Lc X|\Lc^0X}$ (independent on $\eta$) and yielding a map
$$\tau d^{-1}:\Omega^{p,cl}_X\to p_*\Omega^{p-2}_{\Lc X|\Lc^0 X}\leqno(1.4.1)$$
of sheaves on the Zariski topology.
We start with a version of the formal Poincar\'e lemma.

\proclaim{(1.4.2) Lemma}
Let $Z$ be a scheme, locally in the \'etale topology isomorphic to 
$\AA^I=\Spec[x_i;i\in I]$ for some possible infinite set $I$.
Let $W\supset Z$ be a formal scheme locally isomorphic to
$\Spf\Oc_Z[[y_1,...y_m]]$.
Let $\alpha\in\Omega^{p-1}_{W|Z}$ be a closed $(p-1)$-form on $W$ whose
 restriction
to $Z$ is 0. Then $\alpha=d\beta$ for an unique $\beta\in\Omega^{p-2}_{W|Z}$.
\endproclaim \qed

Notice that $Z=\Lc^0 X$ satisfies the conditions of the lemma, while $\Lc X$
is a union of ind-schemes corresponding to formal schemes $W$ also satisfying 
the conditions.

Let now $\omega$ be a section of $\Omega^{p,cl}_X$ as before.
Then $\tau^p(\omega)$ is a section of $\Omega^{p-1,cl}_{\Lc X|\Lc^0 X}$.
By lemma 1.4.2, $\tau^p(\omega)=d\beta$ for a unique 
$\beta\in\Omega^{j-2}_{\Lc X|\Lc^0X}$. We define
$$\tau d^{-1}(\omega)=\beta.$$
It is clear that if $\omega=d\eta$ for an algebraic $(p-1)$-form $\eta$,
then $\tau d^{-1}(\omega)=\tau(\eta)$.
In this paper we will be interested in the case $p=2$. 
In this case we will call $\tau d^{-1}(\omega)$ the symplectic action 
functional. 
We will also consider the map
$$S:\Omega^{2,cl}_X\to \Oc^\times_{\Lc X|\Lc^0 X},
\quad
S(\omega)=\exp(\tau d^{-1}(\omega)).\leqno (1.4.3)$$
We call $S$ the exponentiated symplectic action map. Note that $S(\omega)$,
being the exponential of the action, is the quantity entering the path
integral in the Feynman interpretation of  quantum field theory.

\subhead{\bf (1.5) The global loop space and factorization}
\endsubhead

Let $C$ be a smooth algebraic curve, 
and $X$ be a smooth algebraic variety over $\CC$. 
For any surjection $J\twoheadrightarrow I$ of
nonempty finite sets we denote by $J_i$ the preimage of $i\in I$.
To such a surjection one associates, in a standard way,
the "diagonal" embedding
$$\Delta^{(J/I)}: C^I\to C^J.$$
Let $U^{(J/I)}\subset C^J$ be the locus of $(c_j)_{j\in J}$
such that $c_j\neq c_{j'}$ 
whenever the images of $j,j'$ in $I$ are different.
We denote by $$j^{(J/I)}:U^{(J/I)}\to C^J$$ the embedding.
For $I=J$ it yields an open embedding 
$$j^{(I/I)}:C^{I,\neq}=U^{(I/I)}\to C^I.$$
For each $I$ we have constructed in [KV1] 
an ind-scheme $\Lc(X)_{C^I}$ 
over $C^I$ with the following properties: 

\vskip .2cm 

\itemitem{(a)} For $|I|=1$, when $C^I=C$, 
the ind-scheme $\Lc(X)_C$ is given by the principal 
bundle construction of Gelfand-Kazhdan : 
$$\Lc(X)_C = \Lc X\times_{K} \widehat{C},$$ 
where $\widehat{C}$ is the bundle of 
formal coordinate systems on $C$. Further, $\Lc(X)_C$ has a
connection along $C$ induced by the action of $d/dt\in\gen$. 

\vskip .2cm 

\itemitem{(b)} The $\Lc(X)_{C^I}$ have a natural structure of a 
factorization monoid in the category of ind-schemes. 
More precisely, they are equipped with flat connections along the
$C^I$. Further,
for every surjection $J\twoheadrightarrow I$
there is

\vskip .2cm 

\itemitem{-} an isomorphism of $C^I$-ind-schemes
$\mu^{(J/I)}:\Delta^{(J/I)*}\Lc(X)_{C^J}\to\Lc(X)_{C^I}$,

\vskip .2cm 

\itemitem{-} an isomorphism of $U^{(J/I)}$-ind-schemes
$$\kappa^{(J/I)}:j^{(J/I)*}\biggl(\prod_{i\in I}\Lc(X)_{C^{J_i}}\biggr)\to
j^{(J/I)*}\Lc(X)_{C^J}.$$

\vskip .2cm 

Further, these data are compatible in a natural way.

\vskip .2cm 

There is also a scheme $\Lc^0(X)_{C^I}\subset\Lc(X)_{C^I}$
satisfying the analogs of the 
properties (a) and (b) with $\Lc^0 X$ instead of $\Lc X$. 
The ind-scheme $\Lc(X)_{C^I}$ is a formal thickening
of $\Lc^0(X)_{C^I}$, so it is given by a sheaf of topological rings
$\Oc_{\Lc(X)_{C^I}}$ on $\Lc^0(X)_{C^I}.$ 

For each local function 
$f$ on $\Lc(X)_{C^J}$ 
we write $f|_{C^I}$ for $\Delta^{(J/I)*}(f)$,
and $f|_{U^{(J/I)}}$ for $j^{(J/I)*}(f)$.
Further we omit the identifications
$\mu^{(J/I)},$ $\kappa^{(J/I)}$.

Let $f$ be a local section of $\Oc_{\Lc(X)}$. 
If $f$ is $K$-invariant, 
then by the principal 
bundle construction it gives rise 
to a local function $f_C\in \Oc_{\Lc(X)_C}$. 
For every surjection $J\twoheadrightarrow I$
define a function
$f_{C^J}^{I,\neq}$ on
$$\Lc(X)_{C^J}^{I,\neq}:=
\Lc(X)_{C^J}|_{\Delta^{(J/I)}(C^{I,\neq})}\simeq
\biggl(\prod_{i\in I}\Lc(X)_C\biggr)\biggr|_{C^{I,\neq}}$$
by the formula
$$f_{C^J}^{I,\neq}\bigl((a_i)_{i\in I}\bigr)
=\sum_{i\in i}f_C(a_i).\leqno(1.5.1)$$ 

\proclaim{(1.5.2) Definition} 
A local function $f$ on $\Lc (X)$ 
is called factorizing, if the following properties
hold 


(a) $f$ is $K$-invariant,

(b) for each $J$ there is a local function
$f_{C^J}$ on $\Lc(X)_{C^J}$ 
restricting to $f_{C^J}^{I,\neq}$ on
$\Lc(X)_{C^J}^{I,\neq}$
for every surjection $J\twoheadrightarrow I$.
\endproclaim 

We denote by 
$\Fc\subset p_*\Oc_{\Lc X}$
the subsheaf formed by factorizing functions vanishing on $\Lc^0(X)$.
We now formulate the main result of this paper.

\proclaim{(1.5.3) Theorem}
The image of $\tau d^{-1}:\Omega^{2,cl}_X\to p_*\Oc_{\Lc X}$
coincides with $\Fc$.
\endproclaim


The proof will occupy the rest of the paper. Here we point out one simple
property.

\proclaim{(1.5.4) Lemma} If a factorizing function is $\gen$-invariant,
then it vanishes on $\Lc^0(X)$.

\endproclaim

\noindent {\sl Proof:} It is enough to assume that $X= \Spec (A)$ is affine.
Recall that $\Lc^0(X)$ is the scheme of infinite jets of maps $\AA^1\to X$.
Therefore $\Lc^0(X) = \Spec (A_\Dc)$ where $A_\Dc$ is the differential envelope
of $A$, i.e., the algebra obtained by adjoining to $A$ all iterated formal
derivatives $\partial^n(a)$ of all elements of $A$ which are subject only
to the Leibniz rule, see [BD], (2.3.2). The derivation $\partial$
on $A_\Dc$ corresponds to the action of $d/dt\in\gen$ on $\Lc^0(X)$. 
It is known that the ring of constants $\Ker(\partial) \subset A_\Dc$
is equal to $\CC$. This means that every $\gen$-invariant function on
$\Lc^0(X)$ is constant. The factorization condition implies that this
constant is equal to 0. \qed

\vskip .3cm

\subhead{\bf (1.6) Radon transforms are factorizing} \endsubhead
Here we prove that $$\Im(\tau d^{-1})\subset\Fc.$$
The conditions (a), (b) of Definition 1.5.2 are clear. It remains to prove (c).
We start by recalling from \cite{KV1} 
the functors represented by $\Lc(X)_{C^J}$.
Let $S$ be a scheme and $a_J=(a_j)_{i\in J}$ be a morphism 
$S\to C^J$, so that $a_j:S\to C$, $j\in J$.
Let $\Gamma_j=\Gamma_{a_j}\subset S\times C$ be the graph of $a_j$ and
$\Gamma=\Gamma_J=\bigcup_j\Gamma_j$ be the union.
This is a Cartier divisor in $S\times C$, 
so it is locally given by one equation.
We denote by $\widehat{\Oc}_\Gamma$ the completion of $\Oc_{S\times C}$ along
$\Gamma$, and by $\Kc_\Gamma$ the localization
$\widehat{\Oc}_\Gamma[t^{-1}]$ where $t$ is a (local)  equation of $\Gamma$.
Let also $\Kc_\Gamma^\sqr\subset\Kc_\Gamma$ be the subsheaf formed by sections
whose restriction to $S_{red}\times C$ lies in $\widehat{\Oc}_{\Gamma_\red}$,
where $\Gamma_{red}=\Gamma\cap(S_{red}\times C)$. 
Then, see \cite{KV1}

\proclaim{(1.6.1) Proposition}
(a) Morphisms $h:S\to\Lc(X)_{C^J}$ are in bijection with systems
$(a_J,\varphi)$ where $a_J:S\to C^J$ and $\varphi$ is a morphism of 
locally ringed spaces
$$(\Gamma,\Kc_\Gamma^\sqr)\to(X,\Oc_X).$$

(b) Similarly, morphisms
$k:S\to\Lc^0(X)_{C^J}$ are in bijection with systems
$(a_J,\psi)$ where $a_J$ is as before and $\psi$ is a morphism of 
locally ringed spaces
$$(\Gamma,\widehat{\Oc}_\Gamma)\to(X,\Oc_X).$$
\endproclaim

Notice that to define a function $f_{C^J}\in\Oc_{\Lc(X)_{C^J}}$ 
is the same as to 
define for each morphism $h:S\to\Lc(X)_{C^J}$ with $S$ a scheme,  
a function $h^*f_{C^J}\in\Oc_S$, in a compatible way. 

Let $\omega\in\Omega^{2,cl}_X$ and 
$f=\tau d^{-1}(\omega)\in\Oc_{\Lc X}$.
Let $h$ be as before, 
so in particular, we have morphisms
$a_j:S\to C$, $j\in J$.
For every surjection $J\twoheadrightarrow I$
let $S^{I\neq}\subset S$ be the subscheme 
given by the conditions $a_j=a_{j'}$ iff the image
of $j,j'$ in $I$ are equal.
The functions $f_{C^J}^{I,\neq}$ 
defined in (1.5) give functions
$h^*f_{C^J}^{I,\neq}\in\Oc_{S^{I,\neq}}$, 
and we need to prove that they are the 
restrictions of a regular function  
$h^*f_{C^J}\in\Oc_S$, which is then defined uniquely.
This can be verified by working in the formal 
neighborhood of each point
$s$ in $S^{I,\neq}$. 
But for this we can replace $X$ 
by the union $U$ of the formal neighborhoods of the points
$$\{x_i\}_{i\in I}=\{\varphi(a_j(s))\}_{j\in I}.$$
Further, in $U$
the form $\omega$ can be represented
as $d\eta$ for a 1-form $\eta$. 
So, in the rest of the proof we will assume that $\omega$ is exact, 
hence $$f=\tau(\eta),\quad\eta\in\Omega^1_X.\leqno(1.6.2)$$
Let $\Omega=\Omega^1_{S\times C/S}$.
The morphism
$$\varphi:(\Gamma,\Kc_\Gamma^\sqr)\to(X,\Oc_X)$$
gives rise to a section 
$$\varphi^*\eta\in H^0(\Gamma,\Kc_\Gamma^\sqr\otimes\Omega).$$
Let $\pi:\Gamma\to S$ be the projection.
We now recall the definition and properties of the relative residue map
$$\Res_S:\pi_*(\Kc_\Gamma\otimes\Omega)\to\Oc_S.\leqno(1.6.3)$$
Let $\Gamma^{(m)}\subset S\times C$ be the $m$-th infinitesimal
neighborhood of $\Gamma$, so
$$\widehat{\Oc}_\Gamma=\pro_m\Oc_{\Gamma^{(m)}}\quad\roman{and}\quad
\Oc_{\Gamma^{(m)}}=\Oc_{S\times C}/(t^{m+1})$$
where $t$ is a local equation of $\Gamma$ in $S\times C$.
Recall the natural map
$$\ind_m\underline\Ext^1(\Oc_{\Gamma^{(m)}},\Omega)\to\underline H^1_\Gamma(\Omega),
\leqno(1.6.4)$$
see \cite{C, A.2.18}.
Note also that the short exact sequence
$$0\lra\widehat{\Oc}_\Gamma{\buildrel t^{m+1}\over\lra}
\widehat{\Oc}_\Gamma\lra\Oc_{\Gamma^{(m)}}\lra 0$$
yields, after passing to $\Ext(-,\Omega)$, that
$$\Ext^1(\Oc_{\Gamma^{(m)}},\Omega)=
\widehat{\Oc}_\Gamma\otimes\Omega/t^{m+1}(\widehat{\Oc}_\Gamma\otimes\Omega)\simeq
t^{-m-1}(\widehat{\Oc}_\Gamma\otimes\Omega)/\widehat{\Oc}_\Gamma\otimes\Omega,$$
and therefore the LHS of (1.6.4) is identified with
$$\Kc_\Gamma\otimes\Omega/\widehat{\Oc}_\Gamma\otimes\Omega.$$
We get a morphism of sheaves
$$\Kc_\Gamma\otimes\Omega{\buildrel P\over\to}\underline H^1_\Gamma (\Omega).\leqno(1.6.5)$$
We now compose $P$ with the trace map of the Grothendieck duality theory,
see \cite{C, A.2}, [Ha],
$$\tr_{\Gamma/S}:\pi_*\underline H^1_\Gamma(\Omega)\to\Oc_S.$$
Recall that $tr_{\Gamma/S}$ is obtained from the canonical adjunction map
$$\pi_*\pi^!\Oc_S\to\Oc_S$$
via the map
$\underline H^1_\Gamma(\Omega)\to\pi^!\Oc_S.$
So we define the function
$$h^*f_{C^J}=\tr_{\Gamma/S}(P(\varphi^*\eta))\in\Oc_S.$$
It is clear that these functions are 
compatible for different $h:S\to\Lc(X)_{C^J}$,
and thus define a function
$f_{C^J}\in\Oc_{\Lc(X)_{C^J}}.$
We now verify that $f_{C^J}$ restricts to
$f^{I,\neq}_{C^J}$ on
$\Lc(X)_{C^J}^{I,\neq}$.
Working with $h$ as before,
it suffices to identify $h^*f_{C^J}$ on $S^{I,\neq}$.

We can replace $S$ by $S^{I,\neq}$ and assume that 
$a_{j}=a_{j'}=a_i$ for each $j,j'\in J$ 
whose images in $I$ are equal to $i$.
Further the graphs $\Gamma_i$ of $a_i$, for $i\in I$
are disjoint. 
Representing locally $\varphi^*\eta$ as $\eta_0/t^d$ with 
$\eta_0\in\widehat{\Oc}_\Gamma\otimes\Omega$, we find that,
 see \cite{C, A.2.1},
$$tr_{\Gamma/S}(P(\varphi^*(\eta)))=
\bigl[{\eta_0\over t^d}\bigr].\leqno(1.6.6)$$
In virtue of \cite{C, A.1.5} 
the "residue symbol" on the RHS of (1.6.6)
is nothing but the sum of the residues
$$\sum_{i\in I}\Res_{a_i}(\varphi^*\eta).$$
This is precisely the definition of the function 
$f_{C^J}^{I,\neq}$ in (1.5.1), (1.6.2).

\newpage

\centerline{\bf 2. The sheafified Heisenberg module.} 

\vskip 1cm

\subhead{\bf (2.1) Reminder on vertex algebras}\endsubhead

In this paper we will consider only $\ZZ_{\ge 0}$-graded vertex algebras as in 
\cite{GMS2}.
 Such an algebra is a $\ZZ_{\ge 0}$-graded $\CC$-vector space
$V=\bigoplus_{n\ge 0}V_n$ equipped with a distinguished element $1\in V_0$,
an endomorphism $\partial:V\to V$ of degree one, and a family of
 bilinear operations
$(x,y)\mapsto x\circ_n y$, such that 
$V_i\circ_n V_j\subset V_{i+j-n-1}$,
subject to the axioms in \cite{GMS2, Def 0.4}.
One defines $Y(x,t)$ to be the operator formal series
$$Y(x,t):y\mapsto\sum_n(x\circ_ny)t^{-n-1}.$$
Recall that a morphism of $\ZZ_{\ge 0}$-graded vertex algebras 
$T:V\to W$ is a homogeneous linear map such that 
$$T(Y(x,t)y)=Y(Tx,t)Ty,
\quad
T(1)=1,
\quad\forall x,y\in V.$$

\noindent{\bf (2.1.1) Examples.}
(a) If $A$ is a commutative algebra with a derivation $\partial$, 
it is made into a vertex algebra with
$$Y(a,t)b=\exp(t\,\partial)(a)\cdot b.$$
Such vertex algebras will be called commutative (or holomorphic).
They are characterized by the property that
$Y(a,t)\in\End(A)[[t]]$.
Note that in this case $a\circ_{-1}b=a\cdot b$ is the usual product in $A$.

\vskip2mm

(b) Let $N\ge 1$ and $A_N$ be the Heisenberg algebra with generators 
$a^i_n$, $b^i_n$, $i=1,...,N$, $n\in\ZZ$, subject to
$$[a^j_m,b^i_n]=\delta_{i,j}\delta_{m,-n},$$
all other brackets being zero.
Let $V_N$ be the cyclic $A_N$-module generated by a vector $\1b$ subject to
$$b^i_{<0}\cdot\1b=a^i_{\le 0}\cdot\1b=0.\leqno(2.1.2)$$
It is known that $V_N$ has a structure of a graded vertex algebra
called the Heisenberg vertex algebras such that
$$Y(a_1^i{\bold 1},t)x=\sum_na^i_nx\,t^{n-1},\quad
Y(b_0^i{\bold 1},t)x=\sum_nb^i_nx\,t^{n},$$
and
$$\partial(a^i_{n+1}{\bold 1})=(n+1)a^i_{n+2}{\bold 1},\quad 
\partial(b^i_n{\bold 1})=(n+1)b^i_{n+1}{\bold 1},
\quad n\ge 0,$$
see, e.g., \cite{K}, \cite{MSV}.

\vskip3mm

\noindent{\bf (2.1.3) Definition.}
A filtered vertex algebra is a vertex algebra $V$ equipped with a filtration
$\{F_i V\}_{i\ge 0}$ such that
$$(F_iV)\circ_n(F_jV)\subset F_{i+j-n-1}V.$$
In this case $gr^FV=\bigoplus_iF_iV/F_{i-1}V$
inherits the structure of a vertex algebra.

\vskip3mm

\noindent{\bf (2.1.4) Example.}
Let us introduce a filtration on the vertex algebra $V_N$ from 4.1.1.(2).
A basis of $V_N$ is formed by the monomials
$$a^{i_1}_{m_1}\cdots a^{i_p}_{m_p}
b^{j_1}_{n_1}\cdots b^{j_q}_{n_q}\cdot\1b,\quad
m_i>0,\,n_j\ge 0.$$
We define $F_lV_N$ to be the span of the monomials above with $p\le l$.
Note that 
$$F_0V_N=\CC[b^j_n;j=1,...N,n\ge 0]$$
is in fact a holomorphic vertex algebra. 
It corresponds to the derivation of the commutative algebra
$\CC[b_n^j]$ such that $\partial(b^j_n)=(n+1)b^j_{n+1}$.
Further, $gr^FV_N$ is again holomorphic and identified 
with the commutative algebra
$$\CC[a^i_m,b^j_n;m>0,n\ge 0]$$
equipped with the derivation
$$\partial(a^i_{m})=ma^i_{m+1},\quad 
\partial(b^i_n)=(n+1)b^i_{n+1}.$$

\subhead{(2.2)  Chiral differential operators on $\AA^N$} \endsubhead

Let $\Oc_{\AA^N}^{\ch}$ be the sheaf of vertex algebras
defined in \cite{MSV}. 
Recall that $b^1,b^2,...b^N$ are the coordinates in $\AA^N$.
Identifying $b^i$ with $b^i_0$,
localization of the $\CC[b_0^1,...b_0^N]$-module $V_N$
yields a quasicoherent Zariski sheaf on $\AA^N$.
See \cite{MSV} for the localization of the vertex structure.
The filtration $F$ on $V_N$ induces a filtration (also denoted $F$) on
$\Oc^{\ch}_{\AA^N}$.
We call $\Oc_{\AA^N}^{\ch}$ the sheaf 
of chiral differential operators on $\AA^N$.

Notice further that $V_N$ is a module over the algebra
$$\CC\bigl[b_n^i;i=1,...N, n\ge 0\bigr]
\bigl[\bigl[b_n^i;i=1,...N, n<0\bigr]\bigr]=
H^0(\Lc^0\AA^N,\Oc_{\Lc\AA^N}).$$
This structure comes from the $A_N$-module structure, because sufficiently high
 degree monomials in the $b^i_{<0}$ vanish on any given vector.
Therefore by localizing we get a sheaf $\Vc_N$ of discrete modules over the
sheaf of topological rings $\Oc_{\Lc\AA^N}$ on the space $\Lc^0\AA^N$.
We have
$$\Oc^{\ch}_{\AA^N}=p_*(\Vc_N),\quad
V_N=H^0(\Vc_N).\leqno(2.2.1)$$
In particular, $\Oc^{\ch}_{\AA^N}$ is a
$p_*\Oc_{\Lc\AA^N}$-module.

\subhead{(2.3) The sheaf of automorphisms of $\Oc^{\ch}_{\AA^N}$} \endsubhead
 
We denote by $\Ac$ the Zariski sheaf on $\AA^N$ formed by automorphisms of
$\Oc^{\ch}_{\AA^N}$ as a sheaf of filtered vertex
algebras that induce the identity on
$\gr^F\Oc^{\ch}_{\AA^N}$.
The following result is due to \cite{GMS1}.

\proclaim{(2.3.1) Proposition}
The sheaf $\Ac$ over $\AA^N$ is identified with $\Omega^{2,cl}_{\AA^N}$.
\endproclaim

For future reference we recall the construction of \cite{GMS1}.
Let
$$\omega=\sum_{i,j}h_{ij}(b^1,...,b^N)db^i db^j$$
be a closed 2-form regular in a Zariski open set in $\AA^N$.
Let us define for convenience the generating functions
$$b^i(t)=\sum_{n\in\ZZ}b_n^it^n,\quad
a^i(t)=\sum_{n\in\ZZ}a_n^it^{n-1},$$
see \cite{MSV, (1.11), (1.17)}.
The formulas (4.4a),(4.4b) of \cite{GMS1} 
(taken for the case $g=\roman{Id}$) yields
the following expression for the automorphism
$T_\omega: \Oc^{\ch}_{\AA^N}\to\Oc^{\ch}_{\AA^N}$
corresponding to $\omega$ by 2.3.1 

$$b^i(t)\mapsto b^i(t),\leqno(2.3.2)$$
$$a^i(t)\mapsto a^i(t)+\sum_{k=1}^Nb^k(t)'h_{ki}(b^1(t),...b^N(t)).\leqno(2.3.3)$$

\noindent{\sl Proof :}
It was proved in loc.\ cit.\ that 2.3.2-3 indeed define an automorphism 
$T_\omega$ of sheaves of graded vertex algebras. 
It is clear that $T_\omega$ preserves the 
filtration $F$ introduced in Example 2.1.4 
and induces the identity on the quotients
of this filtration. This is because $F$ is defined in terms of the number of 
monomials in the $a^i_n$. Thus $T_\omega$ is a section of $\Ac$.

Further, the filtration induced by $F$ on the degree one part 
$$(\Oc^{\ch}_{\AA^N})_1\subset\Oc^{\ch}_{\AA^N}$$
coincides with the two-step filtration introduced in \cite{GMS1, 2.2}.
It now follows from the results of \cite{GMS1} that each automorphism of
$\Oc^{\ch}_{\AA^N}$ preserving the grading and the filtration $F$,
preserves, in particular, the two-step filtration on 
$(\Oc^{\ch}_{\AA^N})_1$ and so is of the form $T_\omega$, 
$\omega\in\Omega^{2,cl}_{\AA^N}$. 
\qed

\proclaim{(2.3.4) Proposition} Let $\omega$ be a local section of
$\Omega^{2,cl}_{\AA^N}$ and 
$S(\omega)=\exp(\tau d^{-1}(\omega))\in\Oc^\times_{\Lc\AA^N}$ 
be the corresponding symplectic action functional. Then
$T_\omega: \Oc^{\ch}_{\AA^N}\to\Oc^{\ch}_{\AA^N}$
becomes equal, after the identification 
$\Oc^{\ch}_{\AA^N}=p_*\Vc_N$,
to the operator of multiplication by $S(\omega)$ on the
$\Oc_{\Lc\AA^N}$-module $\Vc_N$.
\endproclaim

\noindent{\sl Proof :}
Recall that the algebra $H^0(\Oc_{\Lc(\AA^N)})$ 
acts on $$V_N=H^0(\AA^N,\Oc^{\ch}_{\AA^N})=H^0(\Lc(\AA^N),\Vc_N).$$
We can think of the element $\bold{1}\in V_N$ as global section of
the sheaf $\Vc_N$ over $\Lc(\AA^N)$.
Further, relations (2.1.2) means that 
the sheaf of $\Oc_{\Lc(\AA^N)}$-modules $\Vc_N$
is isomorphic to the sheaf of distributions on $\Lc(\AA^N)$
supported on $\Lc^0(\AA^N)$, with
$\1b$ identified with the delta function $\delta=\delta_{\Lc^0\AA^N}$ 
(we'll see later, in (3.2), that this isomorphism is compatible with the
vertex structure).
So the element $a_n^i{\bold 1}\in V_N$ corresponds to the 
derivative 
$\partial\delta/\partial b_{-n}^i$. 
Let now $\Phi = \Phi(b)$ be any element of $H^0(\Oc_{\Lc(\AA^N)})$. 
Then we have 
$$\Phi \cdot {\partial \delta\over\partial b_{-n}^i} = 
\bigl( \Phi|_{\Lc^0\AA^N}\bigr)\cdot 
\biggl({\partial\delta\over\partial b_{-n}^i}\biggr) + 
{\partial \Phi\over \partial b_{-n}^i}\cdot 
\delta.$$ 
Let us now specialize to $\Phi=S(\omega)$.
In this case $\Phi$ is equal to 1 on $\Lc^0\AA^N$ so we get 
$$\Phi\cdot(a_n^i\1b)=
a_{n}^i\1b+{\partial\Phi\over\partial b_{-n}^i}\cdot\delta.$$ 
Next, since $S(\omega) = \exp(\tau(d^{-1}(\omega)))$, we have 
$${\partial \Phi\over\partial b_{-n}^i}=
{\partial\over\partial b_{-n}^i}\tau(d^{-1}(\omega))\cdot\Phi 
=\biggl\langle\tau(\omega),
{\partial\over\partial b_{-n}^i}\biggr\rangle\cdot\Phi.
\leqno (2.3.5)$$ 
Here the expression in the angle brackets is the contraction of the 1-form 
$\tau(\omega)$ 
and the vector field $\partial/\partial b_{-n}^i$ on $\Lc\AA^N$. 

Recall now the 
definition of the transgression map $\tau$ on 2-forms. As before, assume that 
$\omega = \sum h_{ij} db^i db^j$. Let $(b^i_n)$ be a point of $\Lc\AA^N$ 
(with values in some ring) and $(\delta b^i_n)$ be a tangent vector to 
$\Lc\AA^N$ at $(b^i_n)$. We write the generating series 
$$b^i(t) = \sum_n b^i_n t^n,\quad\delta b^i(t) = \sum_n (\delta b^i_n) t^n.$$ 
Then by definition of $\tau$ (Example 1.3.8) we have 
$$\tau(\omega)((b^i_n), (\delta b_n^i)) = \Res\sum_{i,j} h_{ij} 
(b^1(t), ..., b^N(t)) \delta b^j(t) \cdot b^i(t)'dt,$$ 
see (1.3).
We apply this to the RHS of (2.3.5). 
Contracting with $\partial/\partial b^i_{-n}$ means that we take 
$$\delta b^j_m= 
\cases 1,\roman{if} j=i,\, m=-n \cr 
0,\, \roman{otherwise} 
\endcases.$$ 
On the level of generating series this entails 
$$\delta b^j(t)=\cases t^n,\, j=i,\cr 0,\, j\neq i\endcases.$$ 
Substituting this into (2.3.5), we get 
that $\bigl\langle \tau(\omega),\partial/\partial b^i_{-n}\bigr\rangle$ 
is equal to the 
coefficient at $t^{n}$ in $\sum_k h_{kj}(b(t))\cdot b^k(t)'$. 
So
$$\Phi\cdot a^i(t){\bold 1}=
a^i(t){\bold 1}+\sum_{k=1}^Nb^k(t)'h_{ki}(b^1(t),...b^N(t)){\bold 1}=
T_\omega(a_i(t){\bold 1})$$
because $T_\omega$ is a morphism of vertex algebras satisfying (2.3.3).

Next, we verify the equality on elements of the form $b^i_n{\bold 1}$. 
We get, for any $\Phi$ 
$$\Phi\cdot b^i_n{\bold 1}=\Phi b^i_n \delta_{\Lc^0 \AA^N}= b^i_n \Phi 
\delta_{\Lc^0\AA^N} =  b^i_n \bigl( \Phi|_{\Lc^0 \AA^N}\bigr)\cdot 
\delta_{\Lc^0\AA^N}.$$ 
Let us now specialize to $\Phi=S(\omega)$.
Since $\Phi$ is identically equal to 1 on $\Lc^0\AA^N$, the 
above element is equal to $b^i_n{\bold 1}$. 
So
$$\Phi\cdot b^i(t){\bold 1} = b^i(t){\bold 1}=T_\omega(b_i(t){\bold 1}),$$
because $T_\omega$ is a morphism of vertex algebras satisfying (2.3.2).

Finally, we extend the equality from the generators 
$a^i_n{\bold 1}, b^i_n{\bold 1}$ 
of $V_N$ to the entire vertex algebra $V_N$. 
For this, we notice that both transformations we 
consider are in fact automorphisms of vertex algebras, 
and that $V_N$ is generated by
$a^i_n{\bold 1}, b^i_n{\bold 1}$.
More precisely, the transformations $T_\omega$ 
of \cite{GMS1} are proved in op.\ cit.\ to be automorphisms of vertex algebras. 
On the other hand the multiplication by 
$S(\omega) = \exp(\tau(d^{-1}\omega))$ is a morphism of vertex algebras,
because so is multiplication by $\exp(f)$ for any factorizing function $f$,
see (3.3) below, and Radon transforms are factorizing by (1.6).
Alternatively, up to localizing we may assume that
$\omega=d\eta$ for some $\eta\in\Omega^1_{\AA^N}$.
Set $\eta=\sum_i h_idb^i$.
Then $S(\omega)=\exp(\tau\eta)$, and
$$\tau(\eta)(b^i_n)=\Res\sum_ih_i 
(b^1(t), ..., b^N(t))b^i(t)'dt=
\roman{Res}\,Y\biggl(\sum_ih_i(b^1_0, ..., b^N_0)b^i_1{\bold 1},t\biggr)
\leqno (2.3.6)$$ 
(no normal ordering because the $b^i_n$'s commute with each other).
Thus multiplication by 
$S(\omega)$ is a morphism of vertex algebras,
as an instance of the automorphism 
$$\exp(\roman{Res}\, Y(x,t))$$ 
valid for any vertex algebra $V$ and any element $x\in V$
such that the exponential is well-defined,
see \cite{MSV, sec.~1.8}.
The proposition is proved. \qed 

\proclaim{(2.3.7) Corollary}
If $f\in\Oc^\times_{\Lc \AA^N}$ is a local invertible function such that 
the multiplication by $f$ is an automorphism of $\Oc^{\ch}_{\AA^N}$ as a filtered 
graded vertex algebra, then $f=S(\omega)$ for some $\omega\in\Omega^{2,cl}_{\AA^N}$.
\endproclaim

\noindent{\sl Proof:}
By Prop.~2.3.4 there is an $\omega\in\Omega^{2,cl}_{\AA^N}$ such that 
$f-S(\omega)$ acts by 0 on $\Oc^{\ch}_{\AA^N}$.
Since $\Oc_{\Lc^0(\AA^N)}\subset \Oc^{\ch}_{\AA^N}$ is an $\Oc_{\Lc\AA^N}$-submodule
(generated by the vacuum vector ${\bold 1}$) and the automorphism given by $f$
preserves $\bold 1$, we find that $f=1$ on $\Lc^0(\AA^N)$.
 Since both $f$ and $S(\omega)$ are equal to 1 on $\Lc^0\AA^N$,
$f-S(\omega)$ vanishes on $\Lc^0\AA^N$.
Further, $\Oc^{ch}_{\AA^N}$ can be seen as the module of distributions
on $\Lc\AA^N$ supported on $\Lc^0\AA^N$, so vanishing of the operator of 
multiplication by $f-S(\omega)$ implies that all the iterated normal derivatives
of $f-S(\omega)$ vanish, so $f=S(\omega)$.
\qed

\vskip .3cm

\subhead{(2.4) Chiral differential operators in a coordinate chart}
\endsubhead
Let $X$ be a smooth algebraic variety equipped with an \'etale map
$\phi: X\to\AA^N$. We have then an \'etale morphism of schemes
$\Lc^0\phi: \Lc^0(X)\to\Lc^0(\AA^N)$ and a morphism of
ind-schemes $\Lc\phi: \Lc(X)\to\Lc(\AA^N)$. We have a sheaf of
vertex algebras on $X$
$$\Oc^{\ch}_{X,\phi} = \phi^* \Oc^{ch}_{\AA^N} = \bigl(\phi^{-1}\Oc^{\ch}_{\AA^N})\otimes_{
\phi^{-1}\Oc_{\AA^N}}\Oc_X,\leqno (2.4.1)$$
called the sheaf of chiral differential operators (CDO)
 corresponding to the \'etale coordinate chart $\phi$. 
It inherits a grading and a filtration from $\Oc^{ch}_{\AA^N}$. 
This sheaf is also a sheaf of $p_*\Oc_{\Lc(X)}$-modules so
$$\Oc^{\ch}_{X,\phi} = p_* \Vc_{X,\phi}, \quad \Vc_{X,\phi} = (\Lc\phi)^*\Vc_N.\leqno (2.4.2)$$

\proclaim{(2.4.3) Proposition}
If $f\in\Oc^\times_{\Lc X}$ is a local invertible function such that 
the multiplication by $f$ is an automorphism of $\Oc^{\ch}_{X,\phi}$ 
as a filtered graded vertex algebra, 
then $f=S(\omega)$ for some $\omega\in\Omega^{2,cl}_{X}$.
\endproclaim

\noindent {\sl Proof:} As in (2.3.7), this follows from the fact that the sheaf of
automorphisms of $\Oc^{\ch}_{X,\phi}$ as a filtered graded vertex algebra,
is identified with $\Omega^{2, cl}_X$, see [GMS1]. The arguments
in (2.3.1-4) can be repeated verbatim using the \'etale coordinates on $X$
given by $\phi$ instead of the standard coordinates in $\AA^N$. 
\qed

\newpage

\centerline{\bf 3. The factorization structure.} 

\vskip 1cm

\subhead{\bf (3.1) Reminder on factorization algebras and vertex algebras}
\endsubhead
We recall the concept of a factorization algebra, 
see \cite{BD} for more details.

\proclaim{(3.1.1) Definition}
Let $\Ec$ be a quasi-coherent sheaf on $C$.
A  structure of a
factorization algebra on $\Ec$ is a collection
of quasi-coherent $\Oc_{C^I}$-modules $\Ec_I$
for each non-empty finite set $I$,
such that $\Ec_I$ is flat along the diagonal strata,
$\Ec_{\{1\}}=\Ec$, and

(a)  an isomorphism of $\Oc_{C^I}$-modules
$\nu^{(J/I)}\,:\,\Delta^{(J/I)*}\Ec_J\simto\Ec_I$
for every $J\twoheadrightarrow I$, compatible with the compositions
of $J\twoheadrightarrow I$,

(b) an isomorphism of $\Oc_{U^{(J/I)}}$-modules
$$\varkappa^{(J/I)}\,:\,j^{(J/I)*}(\boxtimes_I\Ec_{J_i})
\simto j^{(J/I)*}\Ec_{J}$$
for every $J\twoheadrightarrow I$, compatible with the compositions of
$J\twoheadrightarrow I$ and with $\nu$,

(c) a global section $1_\Ec\in H^0(C,\Ec)$ such that for every
$f\in\Ec$ one has
$1_\Ec\boxtimes f\in\Ec_{\{1,2\}}\subset j_*j^*(\Ec\boxtimes\Ec)$
and $\Delta^*(1_\Ec\boxtimes f)=f$.

\endproclaim

A  morphism of factorization algebras $\{\Ec_I\}_I\to\{\Fc_I\}_I$
is defined in the obvious way :
it is a collection of morphisms of sheaves of $\Oc_{C^I}$-modules 
$f_I:\Ec_I\to\Fc_I$
which is compatible with $\nu$, $\varkappa$, and is such that $f=f_{\{1\}}$
sends $1_\Ec$ to $1_\Fc$.
One of the main results of \cite{BD} is as follows.

\proclaim{(3.1.2) Proposition}
Let $\Ec$ be a quasicoherent sheaf on $C$ 
equipped with a structure of a factorization algebra. 
Let $0\in C$ be a $\CC$-point, and $t$ be a formal coordinate
on $C$ near 0. Then $\Ec_0=\Ec\otimes_{\Oc_C}\Oc_0$ has a natural
structure of a vertex algebra (depending on the choice of $t$).
Further, a morphism of factorization algebras 
yields a morphism of vertex algebras.
\endproclaim

\subhead{\bf (3.2) The factorization 
algebra corresponding to $\Oc^{\ch}_{\AA^N}$}
\endsubhead

We now assume that the curve $C$ is $\AA^1$ with coordinate $t$ and $X=\AA^N$.
We identify $\Lc(X)$, $\Lc^0(X)$ 
with the fibers of $\Lc(X)_C$, $\Lc^0(X)_C$ over $0\in C$.
For each $I$, let
$$q_I:\Lc^0(X)_{C^I}\to X\times C^I$$
be the structural morphism.
The aim of this section is to describe a sheaf of
$\Oc_{\Lc(X)_{C^I}}$-modules $\Vc_{N,I}$ over
$\Lc^0(X)_{C^I}$ 
such that

\vskip1mm

\itemitem{(a)}
the collection 
$\{(q_I)_*(\Vc_{N,I})\}_I$ 
of sheaves of $\Oc_{C^I}$-modules over $X$
form a sheaf of factorization algebras over $X$,

\vskip1mm

\itemitem{(b)}
the identification
$\Lc(X)\simeq\Lc(X)_{C,0}$
with the fiber of $\Lc(X)_C$ at 0
takes $\Vc_N$ to 
$\Vc_{N,\{1\}}|_{\Lc(X)_{C,0}},$

\vskip1mm

\itemitem{(c)}
$p_*(\Vc_N)$ is isomorphic to $\Oc_{\AA^N}^{ch}$
as a sheaf of vertex algebras. 

\vskip2mm

In order to avoid technical complications in proving (c),
we  give a proof which relies on the results of \cite{KV1}.
To do so, we first recall the main steps in the construction of \cite{KV1}.

First, recall that $X=\AA^N$ and $C=\AA^1$ and that, see \cite{KV1, 3.7.3} : 
$$\Lc^0(\AA^N)_{\AA^I}=
\Spec\CC\bigl[\l_i,b_{n,\nu}^j;
i\in I,n\ge 0,j=1,...N,\nu=1,...|I|\bigr],\leqno(3.2.1)$$
$$\Oc_{\Lc(\AA^N)_{\AA^I}}=
\Oc_{\Lc^0(\AA^N)_{\AA^I}}\bigl[\bigl[b^j_{n,\nu};
n<0,j=1,...N,\nu=1,...|I|\bigr]\bigr].\leqno(3.2.2)$$
More precisely, if $R$ is a commutative algebra and $(\l_i,b^j_{n,\nu})$
is a system of elements of $R$, i.e., an $R$-point of
$\Lc^0(\AA^N)_{\AA^I}$, one associates to it a morphism
$$c_I:\Spec(R)\to\AA^I=\Spec\CC[t_i;i\in I]$$
taking $t_i$ to $\l_i$ and the homomorphism
$$\rho:\CC[b^1,...b^N]\to H^0(\Gamma_I,\widehat{\Oc}_{\Gamma_I}),
\quad
b^j\mapsto
\sum_{n=0}^\infty\sum_{\nu=1}^{|I|}b^j_{n,\nu}t^{\nu-1}
\prod_{i\in I}(t-\l_i)^n,$$
and similarly for $\Lc(\AA^N)_{\AA^I}$.

Put also
$$\Lc^\eps_\ell(\AA^N)_{\AA^I}=
\Spec\biggl(\CC\bigl[\l_i,b_{n,\nu}^j;\ell\ge n\bigr]/
(b_{n,\nu_1}^jb_{n,\nu_2}^j...b_{n,\nu_{1+\eps_n}}^j;n<0)
\biggr)
,\leqno(3.2.3)$$
where
$\eps=(\eps_{-1},\eps_{-2},...)$ with
$\eps_n\ge 0$, $\eps_n=0$ for almost all $n$, and
$\nu_1,...\nu_{1+\eps_n}=1,...|I|$.
Then the ind-scheme
$\Lc(\AA^N)_{\AA^I}$ represents the ind-pro-object
$$``\ind_\eps"``\pro_\ell"\Lc^\eps_\ell(\AA^N)_{\AA^I}.$$ 
Since $\Lc_\ell^0(\AA^N)_{\AA^I}$ is a smooth scheme, 
we may consider its canonical bundle.
It yields an object in a category of 
relative $\Dc_{\Lc(\AA^N)_{\AA^I}/\AA^I}$-modules,
see the proof of \cite{KV1, 5.1.2}.
Let $\Cc\Dc\Rc(\omega_{\AA^N})_{\AA^I}$
be its De Rham complex.
It is a complex of sheaves over $\AA^N$. 
Further, by \cite{KV1, 5.1.2, 5.3.1},
the collection $$\{\Cc\Dc\Rc(\omega_{\AA^N})_{\AA^I}\}_I$$
is a sheaf of factorization algebras over $X$, and the fiber of
$\Cc\Dc\Rc(\omega_{\AA^N})_{\AA^1}$ at the point $0\in\AA^1$
is a complex of sheaves of 
$\ZZ_{\ge 0}$-graded vertex algebras which is isomorphic
to the chiral De Rham complex $\Omega^{\ch}_{\AA^N}$ of \cite{MSV}.

By \cite{KV1, 5.5.5}, the complex
$\Gamma(\AA^N,\Cc\Dc\Rc(\omega_{\AA^N})_{\AA^I})$
is also equipped with a left action of the $\CC[\AA^I]$-algebra 
$\widetilde{CD}_I$ generated by
symbols $a^i_{n,\nu}$, $b^i_{n,\nu}$, $\phi^i_{n,\nu}$, 
$\psi^i_{n,\nu}$, 
$i=1,...N$, $n\in\ZZ$, $\nu=1,...|I|$,
subject to
$$[a^j_{m,\mu},b^i_{n,\nu}]=\delta_{i,j}\delta_{\nu,\mu}\delta_{m,-n},
\quad
[\phi^j_{m,\mu},\psi^i_{n,\nu}]_+=
\delta_{i,j}\delta_{\nu,\mu}\delta_{m,-n},
$$
all other brackets being zero.
Recall that $\AA^I\simeq\Spec\CC[\l_i;i\in I].$
This module is generated by a cyclic vector $\1b_I$ subject to
$$b^i_{<0,\nu}\cdot\1b_I=a^i_{\le 0,\nu}\cdot\1b_I=
\phi^i_{<0,\nu}\cdot\1b_I=\psi^i_{\le 0,\nu}\cdot\1b_I=0.$$
By localizing $\Gamma(\AA^N,\Cc\Dc\Rc(\omega_{\AA^N})_{\AA^I})$
we get a sheaf $\Wc_{N,I}$ of discrete 
$\Oc_{\Lc(\AA^N)_{\AA^I}}$-modules over
$\Lc^0(\AA^N)_{\AA^I}.$

Now, consider the subalgebra
$A_{N,I}\subset\widetilde{CD}_I$ generated by
$\{a^i_{n,\nu}, b^i_{n,\nu}\}$.
Under localization the $\CC[\AA^I]$-submodule 
$$A_{N,I}\1b_I\subset\Gamma(\AA^N,\Cc\Dc\Rc(\omega_{\AA^N})_{\AA^I})$$
yields a subsheaf $\Vc_{N,I}\subset\Wc_{N,I}$ of discrete 
$\Oc_{\Lc(\AA^N)_{\AA^I}}$-modules over
$\Lc^0(\AA^N)_{\AA^I}.$
It may be viewed as the sheaf of distributions on $\Lc(\AA^N)_{\AA^I}$
supported on $\Lc^0(\AA^N)_{\AA^I}$. 
The following is immediate from \cite{KV1, 5.4, 5.5}.

\proclaim{(3.2.4) Lemma} 
The collection of sheaves of $\Oc_{C^I}$-modules $\{(q_I)_*(\Vc_{N,I})\}_I$
over $X$ is a subsheaf of factorization algebras of 
$\{\Cc\Dc\Rc(\omega_{\AA^N})_{\AA^I}\}_I$.
Further, the corresponding subsheaf of vertex algebra of 
$H^0(\Omega^{ch}_{\AA^N})$
is isomorphic to $V_N=H^0(\Oc^{\ch}_{\AA^N})$.
\endproclaim

\subhead{(3.3) The factorization algebra corresponding to $\Oc^{ch}_{X,\phi}$}
\endsubhead
Let $\phi: X\to\AA^N$ be an \'etale map as in (2.4). We have then the
morphism of factorization groupoids $\{\Lc_{C^I}(\phi): \Lc_{C^I}(X)\to
 \Lc_{C^I}(\AA^N)\}$.
Using this, we obtain the factorization algebra corresponding
to $\Oc^{\ch}_{X,\phi}$ as a pullback. More precisely, we define
sheaves 
$\Vc_{X,\phi ,I}\subset\Wc_{X, \phi, I}$ of discrete 
$\Oc_{\Lc(X)_{\AA^I}}$-modules over
$\Lc^0(X)_{\AA^I}$ as
$$\Vc_{X,\phi, I} = (\Lc_{C^I}(\phi))^*\Vc_{N, I},\quad \Wc_{X,\phi, I}=
(\Lc_{C^I}(\phi))^*\Wc_{N, I}.
\leqno (3.3.1)$$
It follows from [KV1] and the definition of $\Oc^{ch}_{X,\phi}$ that
the next proposition holds.

\proclaim{(3.3.2) Proposition}
(a) The collections of sheaves $\{(q_I)_*(\Vc_{X,\phi, I})\}_I$ and  
 $\{(q_I)_*(\Wc_{X,\phi, I})\}_I$ form sheaves of factorization algebras on $X$.

(b) The corresponding sheaves of vertex algebras are $\Oc^{\ch}_{X,\phi}$ and
$\Omega^{\ch}_{X,\phi}$.
\endproclaim

\vskip .3cm

\subhead{\bf (3.4) Factorizing functions are Radon transforms}
\endsubhead
We can now finish the proof of Theorem 1.5.3.
We must prove the inclusion
$$\Fc\subset\Im(\tau d^{-1}).$$
as above, we take  $C=\AA^1$. Further, it is enough to work
Zariski locally on $X$ and therefore we can assume that
$X$ possesses an \'etale map $\phi: X\to\AA^N$. We fix such $\phi$.
Fix also a local factorizing function $f$ on $\Lc(X)$ vanishing on $\Lc^0(X)$,
and local functions $f_{C^I}$ over $\Lc(X)_{C^I}$ as in (1.5.2).
By  the definition of $f_{C^I}$, we have that
 $f_{C^I}$ vanishes on $\Lc^0(X)_{C^I}$.
Thus the section $\exp(f_{C^I})$ 
of $\Oc^\times_{\Lc(X)_{C^I}}$ is well-defined,
and multiplication by $\exp(f_{C^I})$ is a morphism
of sheaf of $\Oc_{C^I}$-modules 
$(q_I)_*(\Vc_{X,\phi,I})$
over $X$.
When $I$ varies, such morphisms are 
obviously compatible with the factorization
structure of 
$(q_I)_*(\Vc_{N,I})$
in (3.2.4). 
Further, the unit maps to itself because $\exp(f_C)=1$ on $\Lc^0(X)_C$. 
Thus we get a factorization algebra automorphism.
By (3.1.2) 
multiplication by $\exp(f)$ is an automorphism of the 
sheaf of vertex algebras
$$\Oc_{X,\phi}^{\ch}=p_*(\Vc_{X,\phi}).$$
So, by (2.4.3),
we have $\exp(f)=S(\omega)$ for some closed two-form
$\omega$ on $X$.
Hence $f=\tau d^{-1}(\omega)$,
yielding the desired inclusion.
Note also that the main result of [MSV] says that the sheaves
$\Omega^{\ch}_{X,\phi} = \Omega^{\ch}_X$ are in fact independent of $\phi$
up to a canonical isomorphism 
 Theorem 1.5.3 is proved.

\vskip .2cm

We now finish the proof of Theorem 0.2.1. To see the equivalence of (i)
 and (ii),
note that $\Oc^{ch}_{X,\phi}$ 
is a $p_*\Oc_{\Lc X}$-submodule of $\Omega^{ch}_{X,\phi}$,
as well as a sheaf of vertex subalgebras. 
So if the operator of multiplication by $\exp(f)$
is an automorphism of $\Omega^{ch}_{X,\phi}$ 
as a sheaf of vertex algebras, it is
an automorphism of $\Oc^{ch}_{X,\phi}$ as a sheaf of vertex algebras and so
$f$ has the form $\tau(d^{-1}(\omega))$ by (2.4.3). On the other hand,
if $f$ is of the form $\tau(d^{-1}(\omega))$, then
the multiplication by $\exp(f)$ is a vertex automorphism by the
same argument as in (2.3.6) (proof of Proposition 2.3.4). 

 The equivalence of (i) and
(iii) follows from Theorem 1.5.3 and the next lemma.

\proclaim {(3.4.1) Lemma} If $f$ is 2-factorizable, then
the multiplication by $\exp(f)$ is an automorphism of $\Oc^{\ch}_{X,\phi}$
and $\Omega^{\ch}_{X,\phi}$. 
\endproclaim

\noindent {\sl Proof:} We need to prove that the multiplication by
$\exp(f)$ preserves the operations $a\circ_n b$, $n\in\ZZ$. 
But  in the language of
factorization algebras 
these operations are obtained by expanding the factorization
isomorphism on $C^2$ near the diagonal. See [BD] (3.5.14). 
So the existence of $f_{C^2}$ implies that the $a\circ_n b$ are
preserved. \qed

\proclaim {(3.4.2) Corollary} 
A 2-factorizable function is factorizable.\qed
\endproclaim

This is natural to expect if we view factorizability
as an infinitesimal analog of additivity in $n$-punctured disks
for all $n\geq 2$.  Any such disk can be decomposed into 2-punctured
ones (pairs of pants) so additivity in pairs of pants implies
additivity in $n$-punctured disks.

\vskip 2cm

\noindent M.K.: Department of Mathematics, Yale University,
10 Hillhouse Avenue, New Haven CT 06520 USA, email:
$<$mikhail.kapranov\@yale.edu$>$.

\vskip .2cm

\noindent E.V.: D\'epartement de Math\'ematiques, Universit\'e
de Cergy-Pontoise, 2 Av. A. Chauvin, 95302 Cergy-Pontoise Cedex, France,
email: $<$ eric.vasserot\@math.u-cergy.fr$>$

\bye